\documentclass{elsart}
 \usepackage{amsmath} 
 \usepackage{amssymb} 

\begin{document}
\begin{frontmatter}
\title{A Basic Elementary Extension of the Duchet-Meyniel Theorem}
\author[SDU]{Anders Sune Pedersen}
\ead{asp@imada.sdu.dk}
\author[SDU]{Bjarne Toft\corauthref{cor}}
\corauth[cor]{Corresponding author.}
\ead{btoft@imada.sdu.dk}

\address[SDU]{Department of Mathematics and Computer Science, University of Southern Denmark, DK-5230 Odense M, Denmark}



\begin{abstract}
  The Conjecture of Hadwiger implies that the Hadwiger number $h$
  times the independence number $\alpha$ of a graph is at least the
  number of vertices $n$ of the graph. In 1982 Duchet and Meyniel
  proved a weak version of the inequality, replacing the independence
  number $\alpha$ by $2\alpha-1$, that is,
$$(2\alpha-1)\cdot h \geq n.$$
In 2005 Kawarabayashi, Plummer and the second author published an
improvement of the theorem, replacing $2\alpha - 1$ by $2\alpha - 3/2$
when $\alpha$ is at least $3$. Since then a further improvement by
Kawarabayashi and Song has been obtained, replacing $2\alpha - 1$ by
$2\alpha - 2$ when $\alpha$ is at least $3$.

In this paper a basic elementary extension of the Theorem of Duchet
and Meyniel is presented. This may be of help to avoid dealing with
basic cases when looking for more substantial improvements. The main
unsolved problem (due to Seymour) is to improve, even just slightly,
the theorem of Duchet and Meyniel in the case when the independence
number $\alpha$ is equal to $2$. The case $\alpha = 2$ of Hadwiger's
Conjecture was first pointed out by Mader as an interesting special
case.
\end{abstract}


\end{frontmatter}
\section{Introduction and Notation}
The Hadwiger number $h(G)$ of a graph $G$ is the maximum $k$ for
which $G$ has the complete graph $K_k$ as a minor. In 1942 Hadwiger
suggested that $h(G) \geq \chi(G)$, where $\chi(G)$ is the chromatic
number of $G$. This conjecture of Hadwiger is still open.

In any colouring of $G$ each colour is used at most $\alpha(G)$
times, where $\alpha(G)$ is the maximum number of independent
vertices in $G$. Hence $\chi(G)\cdot \alpha(G) \geq n$, where $n$
is the number of vertices of $G$. If Hadwiger's Conjecture is true,
then $h(G) \cdot \alpha(G) \geq n$. This weaker form of Hadwiger's
Conjecture is also unresolved. In 1982 the following result was
obtained:

\begin{thm}[The Theorem of Duchet and Meyniel \cite{DuchetMeyniel82}]
Let $G$ be a graph on $n$ vertices with independence number $\alpha$
and Hadwiger number $h$. Then $$(2\alpha-1)\cdot h \geq n.$$
\end{thm}

It was observed by Maffray and Meyniel \cite{MM} that equality holds
in $(2\alpha-1)\cdot h \geq n$ if and only if $\alpha = 1$, i.e. if
and only if $G$ is complete.

Around the same time, but independently from Duchet and Meyniel,
Woodall \cite{W} divided Hadwiger's Conjecture into subconjectures,
one of which is $h(G) \cdot \alpha(G) \geq n$. In this connection
Woodall proved a result very similar to the Theorem of Duchet and
Meyniel:

\begin{thm}[The Theorem of Woodall \cite{W}]
Let $G$ be a graph on $n$ vertices with at least one edge and with
independence number $\alpha$ and Hadwiger number $h$. Then
$$2\alpha\cdot (h-1) \geq n.$$
\end{thm}

The main purpose of this paper is to provide a common generalization
of the Theorem of Duchet and Meyniel and the Theorem of Woodal,
hoping in this way to cover general cases with a best possible
result, thus being able to avoid consideration of these cases in
future investigations.

All graphs considered in this paper are assumed to be simple and
finite. The \emph{clique number} $\omega (G)$ of a graph $G$ is the
cardinality of a maximum clique of $G$, respectively. For any
undefined concepts the reader may refer to~\cite{MR2368647}. Given
some graph $G$ and graph parameter $\mu$ we may, for ease of notation,
write $\mu$ for the value $\mu (G)$ when no confusion is possible.

\section{The Main Theorem}
The following theorem we call the basic elementary extension of the
Duchet-Meyniel Theorem. The term basic refers to its coverage of
general cases, and the term elementary to the fact that the proof
uses only induction and the original idea of Duchet and Meyniel
\cite{DuchetMeyniel82}.

\begin{thm}[The Main Theorem]
Let $G$ be a graph on $n$ vertices with at least one edge and at least
one missing edge (i.e. $G$ is neither edge-empty nor complete). Let
$\omega, \alpha$ and $h$ denote the order of a largest complete
subgraph of $G$, the size of a largest independent set of $G$ and the
Hadwiger number of $G$, respectively. Then
\begin{equation}
(2 \alpha - 1) \cdot (h - 1) + 3 \geq n + \omega,
\label{eq:elem}
\end{equation}
where equality is obtained if and only if
\begin{itemize}
\item[(i)] $G$ is a non-empty forest with a perfect matching
\end{itemize}
 or
\begin{itemize}
\item[(ii)] $G$ contains two disjoint $K_{n/2}$, possibly with some
  edges between them, such that $h = n/2$.
\end{itemize}
\label{th:elementary-Duchet-Meyniel}
\end{thm}

In the inequality of Duchet and Meyniel a $1$ is thus subtracted
from both factors $2\alpha$ and $h$ on the left hand side,
diminishing it by $2\alpha + h - 1$, moreover $3$ is added on the
left hand side, whereas $\omega$ is added to $n$ on the right hand
side. For the first case with equality we have $\omega = h = 2$ and
$\alpha = n/2$. For the second case with equality we have $h =
\omega = n/2$ and $\alpha=2$. 

The idea to replace $n$ by $n+\omega$ was first used in \cite{PST}
and is due to M. Stiebitz.
\begin{cor}
For $\alpha \geq 3$ and $\omega \geq 3$ $$(2\alpha-1)\cdot (h-1)+
2\geq n+\omega.$$
\end{cor}
The proof of Theorem~\ref{th:elementary-Duchet-Meyniel} is by
induction and consists of many cases. Since the theorem does not hold
for complete graphs nor for edge-empty graphs one needs to be careful
when applying the induction hypothesis. This gives rise to the many
cases. There is a considerable overlap between cases in the following
proof, so maybe a shorter, more elegant proof can be designed. The
main part of the proof of Theorem~\ref{th:elementary-Duchet-Meyniel}
is given in Section~\ref{sec:ProofOfMainTheorem}, while
Sections~\ref{sec:Heq2orHeq3}, \ref{sec:AlphaTwo}
and~\ref{sec:OmegaEqTwo} settle the special cases $ h \in \{2,3 \}$,
$\alpha = 2$ and $\omega = 2$, respectively. First of all, we
determine the values of the parameters $\alpha$, $\omega$ and $h$ for
the extremal graphs described in (i) and (ii) of
Theorem~\ref{th:elementary-Duchet-Meyniel}.
\begin{lem}
  A graph $G$ on $n$ vertices has $h = 2$ and $\alpha = n/2$ if
  and only if $G$ is a non-empty forest with a perfect matching.
\label{lm:lemma1}
\end{lem}
\begin{pf}
  Suppose that $G$ is a graph on $n$ vertices with $h = 2$ and $\alpha
  = n/2$. Since $h = 2$ excludes any cycles, it follows that $G$ is a
  non-empty forest and, in particular, $G$ is bipartite. Moreover,
  $\alpha = n/2$ implies that the bipartite graph has partition sets
  $A$ and $B$, where $|A| = |B| = n/2$. If for every set $S \subseteq
  A$, the size of the neighbourhood $N(S)$ is at least $|S|$, then it
  follows from Hall's theorem, that there is a matching of $A$ into
  $B$, and, since $|A|=|B|=n/2$, such a matching is indeed a perfect
  matching.  On the other hand, if there exists some set $S \subseteq
  A$ with $|N(S)| < |S|$, then the set $S \cup B \backslash N(S)$ is
  an independent set of $G$ of size $|S| + |B \backslash N(S)| = |S|
  +n/2 - |N(S)| > n/2$, which contradicts the assumption $\alpha
  = n/2$.

  Conversely, any non-empty forest with an edge has $h = 2$ and is
  bipartite, hence $\alpha \geq n/2$. Since $G$ has a perfect
  matching, $\alpha \leq n/2$. Thus, the reverse implication follows.
\end{pf}

\begin{lem}
  A graph $G$ on $n$ vertices has $h = \omega = n/2$ and $\alpha = 2$
  if and only if $G$ consists of two disjoint $K_{n/2}$, possibly with
  some edges between them, such that $h = n/2 \geq 1$.
\label{lm:lemma2}
\end{lem}
\begin{pf}
  Suppose $G$ is a graph on $n$ vertices with $h = \omega = n/2$ and
  $\alpha = 2$. Let $K_\omega$ denote a complete $n/2$-subgraph of
  $G$. Now we just need to show that $G - K_\omega$ is a complete
  $n/2$-graph. The graph $G$ contains at least two vertices and
  $\omega < n$, since $\alpha \geq 2$. Hence there exists at least one
  vertex $v \in G - K_\omega$.  If $G - K_\omega$ consists of just a
  single vertex, then $n=2$, and $G$ has the desired structure. Hence
  we may assume that $G - K_\omega$ contains at least two vertices.

  Suppose that the vertices of $V(G) \backslash V(K_\omega)$ can be
  partitioned into two non-empty sets $X$ and $Y$ such that there is
  no edge joining a vertex of $X$ to a vertex of $Y$. Both the induced
  graphs $G[X]$ and $G[Y]$ must be complete, since otherwise we could
  find an independent set of cardinality $3$. Let $V_X$ denote the
  vertices $z \in V(K_\omega)$ which are adjacent to every vertex of
  $X$, and let $V_Y := V(K_\omega) \backslash V_X$. Consider some vertex $z \in
  V_Y$. Since $z \in V_Y$, there must be some vertex $x \in X$, which
  is not adjacent to $z$, since otherwise $z \in V_X$. If $z$ is not
  adjacent to some vertex of $Y$, say $y$, then $\{x,y,z \}$ is an
  independent set in $G$, which contradicts $\alpha = 2$. Hence every
  vertex of $V_Y$ must be adjacent to every vertex of $Y$. This shows
  that both $G[V_X \cup X]$ and $G[V_Y \cup Y]$ are complete graphs.
  The complete subgraphs $G[V_X \cup X]$ and $G[V_Y \cup Y]$ are
  disjoint and, together, they span $G$. Furthermore, $\omega = n/2$,
  and so we must have $|V_X \cup X| = n/2$ and $|V_Y \cup Y| = n/2$.
  Hence $G$ has the desired structure.

  Now suppose that $G - K_\omega$ is connected. The graph $G -
  V(K_\omega)$ contains at least two vertices. Suppose that two
  vertices, say $x$ and $y$, of $G - K_\omega$ are non-adjacent.  Any
  vertex $z \in V(K_\omega)$ must be adjacent to $x$ or $y$.  Thus,
  since $G - K_\omega$ is connected, contracting $V(G) - V(K_\omega)$
  in $G$ into a single vertex results in a complete graph on $1 + n/2$
  vertices, which contradicts the assumption $h = n/2$. This shows
  that $G - K_\omega$ must be a complete graph on $n/2$ vertices, and,
  again, $G$ has the desired structure.

  For the reverse implication, if $G$ contains of two disjoint
  complete $n/2$-graphs and $h = n/2 \geq 1$, then $n/2 = h \geq
  \omega \geq n/2$, implying that $\omega = n/2$.  Moreover, $\alpha
  \geq 2$, since $\omega < n$, and $\alpha \leq 2$ from the given
  structure of $G$. Thus, $\alpha = 2$. This proves the reverse
  implication, and so the proof is complete.
\end{pf}

It is an unsolved problem if the graphs characterized in
Lemma~\ref{lm:lemma2} may be recognized by a polynomial-time
algorithm.

\section{Graphs with $h = 2$ or $h=3$}
\label{sec:Heq2orHeq3}
If $G$ is a graph with $h =2$, then $G$ has no cycles and is a forest,
and hence bipartite. Therefore $\alpha \geq n/2$ and
$$(2\alpha - 1)(h-1) + 3 \geq n + 2 = n + \omega,$$
where equality occurs if and only if $\alpha = n/2$, and so it follows from
Lemma~\ref{lm:lemma1} that Theorem~\ref{th:elementary-Duchet-Meyniel}
holds for $h=2$.

Suppose $G$ is a graph with $h = 3$. If $n = 3$, then $G \simeq K_3$
which is excluded by the assumptions of
Theorem~\ref{th:elementary-Duchet-Meyniel}, so $n \geq 4$.  According
to Hadwiger's theorem \cite{H}\footnote{Short proofs of
  Hadwiger's theorem were given by D. R. Woodall in \cite{MR1147805}
  and M. Stiebitz in \cite{MR1411244}.}, any $4$-chromatic graph
contains $K_4$ as a minor, and so, since $h = 3$, it follows that $G$
is $3$-colourable.  Therefore, $\alpha \geq \lceil n/3 \rceil$ and
\begin{eqnarray}
\nonumber (2\alpha - 1)(h - 1) + 3 & \geq & (2 \lceil n/3 \rceil - 1)2 + 3 \\
\nonumber         &   =     & 4 \lceil n/3 \rceil + 1 \\
\nonumber         &   \geq  & \lceil 4n/3 \rceil + 1 \\
\nonumber         &   =     & n + \lceil n/3 \rceil + 1 \\
\label{eq:349563} & \geq    & n + 3.
\end{eqnarray}
Equality occurs in~\eqref{eq:349563} if and only if $\omega = 3$,
$\alpha = \lceil n/3 \rceil = 2$ and $4\lceil n/3 \rceil = \lceil 4n/3
\rceil$. Now $\lceil n/3 \rceil = 2$ implies $n \in \{ 4,5,6 \}$,
while $4\lceil n/3 \rceil = \lceil 4n/3 \rceil$ implies that $n$ is a
multiple of three. Thus, equality occurs in~\eqref{eq:349563} if and
only if $n = 6$, $\omega = 3$ and $\alpha = 2$, in which case it
follows from Lemma~\ref{lm:lemma2}, that $G$ contains two disjoint
complete $3$-graphs. This shows that
Theorem~\ref{th:elementary-Duchet-Meyniel} holds for $h=3$.

\section{Graphs with $\alpha = 2$}
\label{sec:AlphaTwo}
The case $\alpha = 2$ of the main theorem follows from
Lemma~\ref{lm:lemma2} and the following result:
\begin{thm}
  Let $G$ denote a graph on $n$ vertices, and let $\omega, \alpha$ and
  $h$ denote the order of a largest complete subgraph of $G$, the size
  of a largest independent set of $G$ and the Hadwiger number of $G$,
  respectively. Let $\alpha = 2$. Then
\begin{equation}
3 h \geq n + \omega,
\label{eq:alpha2}
\end{equation}
where equality occurs if and only if $h = \omega = n/2$.
\label{th:alpha2}
\end{thm}
\begin{pf}
  Let $K_\omega$ denote a maximum clique of $G$. Since $\alpha = 2$, of course
  $n \geq 2$. If $n = 2$, then $G \simeq \overline{K_2}$, and the
  theorem holds in this case. Suppose that $n \geq 3$. Then, clearly,
  $h \geq 2$, since $\alpha = 2$. If $h = 2$, then $\alpha \geq n/2$,
  and, since $\alpha = 2$, we obtain $n \leq 4$. In this
  case~\eqref{eq:alpha2} holds with equality if and only if
  $h=\omega=n/2$. Considering $n \in \{ 3, 4, 5, 6 \}$ and $h \geq 3$,
  we obtain strict inequality in~\eqref{eq:alpha2} unless $n = 6$ and
  $h = \omega = 3$, in which case we obtain equality, as claimed. If
  $\omega \geq n/2$, then $3h \geq 3 \omega \geq n + \omega$, where
  equality occurs if and only if $h = \omega = n/2$. Thus, we may
  assume $n \geq 7$, $\omega < n/2$ and that the statement of the
  theorem is true for any graph $H$ of order $< n$ and $\alpha
  (H) = 2$.

\begin{itemize}
\item[(1)] Suppose that $G - K$ is disconnected for some complete
  subgraph $K$ of $G$ ($K$ may even denote the empty graph). In this
  case, our assumption that $\alpha = 2$ implies that $G - K$
  must consist of exactly two components and these components must be
  complete graphs, say $K_a$ and $K_b$. For any vertices $x \in
  V(K)$, $y \in V(K_a)$ and $z \in V(K_b)$, the vertex $x$ must
  be adjacent to at least one of $y$ and $z$. Let $A \subseteq
  V(K)$ consist of the vertices $x \in V(K)$, which are
  adjacent to every vertex of $V(K_a)$, and let $B := V(K) -
  A$. Now any vertex $x \in B$ is not adjacent to every vertex of
  $K_a$, say $x$ is not adjacent to $y \in K_a$. Since no set $\{ x,
  y, z \}$, where $z \in V(K_b)$, is an independent set of $G$, it
  follows that $x$ and $z$ are adjacent.  This shows that every vertex
  of $B$ is adjacent to every vertex of $K_b$. Hence $G[A \cup
  V(K_a)]$ and $G[B \cup V(K_b)]$ are disjoint complete graphs, and, since
  they cover $G$, at least one of them must contain at least $\lceil
  n/2 \rceil$ vertices, and so $\omega \geq n/2$, which contradicts
  our assumption $\omega < n/2$.

\item[(2)] Suppose that $G - K_\omega$ is connected. The graph $G -
  K_\omega$ cannot be a complete graph, since, by the assumption
  $\omega < n/2$, this would imply $n = \omega + n(G - K_\omega) \leq
  2 \omega < n$. Since $G - K_\omega$ is connected, $G - K_\omega$
  must contain an induced $3$-path, say $P_3 : xyz$. Moreover, $G -
  K_\omega$ has at least $\lceil n/2 \rceil \geq 4$ vertices, since
  $\omega < n/2$ and $n \geq 7$. The assumption $\alpha = 2$ implies
  that any vertex of $V(G) \backslash \{x,z\}$ is adjacent to $x$ or
  $z$; we say that $P_3$ \emph{dominates} $G$.  Thus, any complete
  order $k$ minor of $G - P_3$ can, by contracting the two edges of
  the $3$-path $P_3$, be extended to a complete order $k+1$ minor of
  $G$. Define $G' := G - P_3$ and let $\alpha' := \alpha (G')$, $n' :=
  n(G')$, $h' = h(G')$ and $\omega' := \omega(G')$.  Thus, we have $h'
  + 1 \leq h$ and $\omega' = \omega$.  The graph $G'$ cannot be
  complete, since it has more than $\omega$ vertices. Hence we may
  apply the induction hypothesis to $G'$, and obtain
\begin{equation}
3h \geq 3(h' + 1) = 3h' + 3 \geq n' + \omega' + 3 = n + \omega,
\label{eq:2345620123}
\end{equation}
which is the desired inequality. Equality in~\eqref{eq:2345620123}
implies $h=h'+1$ and, by induction, $h' = \omega' = n' /2$, and so
$\omega = (n-3)/2$ and $h = \omega + 1 = (n-1)/2$. Moreover, $G' := G
- P_3$ is of the exceptional type described in Lemma~\ref{lm:lemma2},
in particular, it contains two disjoint complete $\omega$-graphs; let
$V_1$ and $V_2$ denote the vertices of those two complete
$\omega$-graphs.

If each vertex of $V_1$ has a neighbour in $V_2$, then, by contacting
the vertices of $V_2$ and the vertices of $V(P_3)$ into two distinct
vertices, a complete $(\omega + 2)$-minor of $G$ is obtained, which
contradicts the assumption $h = \omega + 1$. This shows that, in $G$,
at least one vertex of $V_1$, say $x'$, has no neighbour in $V_2$.
Similarly, we may assume that some vertex $z' \in V_2$ has no
neighbours in $V_2$.

Since both $G[V_1]$ and $G[V_2]$ are maximum cliques in $G$, none of
the vertices $x,y$ and $z$ are adjacent to every vertex of $V_1$ or
$V_2$. Let $y_1 \in V_1$ and $y_2 \in V_2$ denote non-neighbours of
$y$.

Now, since $\alpha = 2$, the set $\{ y, y_1, z' \}$ is not
independent, and so $yz' \in E(G)$. Similar arguments show that each
vertex of $\{x,y,z \}$ is adjacent to both $x'$ and $z'$. In
particular, $G[y,z,z'] \simeq K_3$ and $\omega \geq 3$.

If some vertex $t \in V_1 \backslash \{ x' \}$ is adjacent to $y$,
then sets $\{y,t, y_1 \}$ and $\{ x, x', z \}$ both induce dominating
$3$-paths, and both are disjoint from $V_2$. Thus, by contracting
$\{y,t, y_1 \}$ and $\{ x, x', z \}$ into two distinct vertices, a
$(\omega + 2)$-minor is obtained, a contradiction. A similar argument
shows that no vertex $t \in V_2 \backslash \{ z' \}$ is adjacent to
$y$. Thus, we obtain $N_G(y) = \{ x,x',z,z' \}$. The graph $G'' := G -
x - z - z'$ has at least two non-adjacent vertices, and so the
induction hypothesis applies to $G''$. Since we are assuming $G$ to be
a graph for which equality is obtained in~\eqref{th:alpha2}, it
follows, exactly as in~\eqref{eq:2345620123}, that $G''$ is a graph for
which equality is obtained in~\eqref{th:alpha2}, and so, by induction,
$G''$ contains two disjoint complete $\omega$-graphs. In particular,
$y$ has at least $\omega-1$ neighbours in $G''$. Thus, in $G$, the
vertex $y$ is adjacent to at least three vertices $x, z, z' \in V(G)
\backslash V(G'')$ and $\omega - 1$ vertices of $V(G'')$, and so $4 =
\deg_G (y) \geq (\omega -1) + 3$, which implies $\omega \leq 2$,
contradicting the fact that $\omega \geq 3$. This completes the proof.
\end{itemize}
\end{pf}

\section{Graphs with $\omega  = 2$}
\label{sec:OmegaEqTwo}
The case $\omega = 2$ of the main theorem follows from
Lemma~\ref{lm:lemma1} and the following result:
\begin{thm}
  Let $G$ be a graph on $n$ vertices, and let $\omega$, $\alpha$ and
  $h$ denote the order of a largest complete subgraph of $G$, the size
  of a largest independent set of $G$ and the Hadwiger number of $G$,
  respectively. Let $\omega = 2$. Then
\begin{equation}
(2 \alpha - 1)(h - 1) + 1 \geq n,
\label{eq:omega2}
\end{equation}
where equality occurs if and only if $h=2$ and $\alpha = n/2$.
\label{th:omega2}
\end{thm}
\begin{pf}
  Let $G$ denote a graph of order $n$ with at least one edge. If
  $\alpha = n/2$ and $h = \omega = 2$, then equality holds
  in~\eqref{eq:omega2}.

  Now for the converse. If $G$ is a complete graph, then $G \simeq
  K_2$, $h=2$, $\alpha =1$ and we have equality in~\eqref{eq:omega2}.
  Thus, we may assume $\alpha \geq 2$. It follows from
  Section~\ref{sec:Heq2orHeq3} that the desired statement holds in the
  cases $h \in \{2,3\}$. Thus, we assume $h \geq 4$ and proceed by
  induction on the order of the graph $G$ to obtain strict inequality
  in~\eqref{th:omega2}.

  Suppose that $G$ is disconnected, and let $G_1$ and $G_2$ be two
  non-empty subgraphs of $G$ such that $G_1 \cup G_2 = G$ and $V(G_1)
  \cap V(G_2) = \emptyset$. Let $\alpha_i := \alpha (G_i)$, $h_i :=
  h(G_i)$ and $n_i := n(G_i)$. Then $\alpha = \alpha_1 + \alpha_2$, $n
  = n_1 + n_2$ and $h =  \max \{ h_1, h_2 \}$. We may assume that
  $h_1 = h \geq 4$, and so the induction hypothesis is applicable to
  $G_1$:
\begin{equation}
(2 \alpha_1 - 1)(h_1 - 1) + 1 \geq n_1.
\label{eq:2938545}
\end{equation}
Now there are two cases to consider depending on whether $G_2$ is
edge-empty or not. Firstly, if $G_2$ is edge-empty, then $\alpha_2 =
n_2$ and, as the following computations show, we obtain
strict inequality in~\eqref{eq:omega2}.
\begin{eqnarray*}
(2\alpha - 1)(h - 1) + 1 & = & (2\alpha_1 + 2n_2 - 1)(h - 1) + 1 \\
& \geq & (2\alpha_1 - 1)(h_1 - 1) + 1 + 2n_2(h - 1) \\
& \geq & n_1 + 6 n_2 > n.
\end{eqnarray*}
Secondly, if $G_2$ is not edge-empty, then the induction hypothesis
also applies to $G_2$, and we obtain a bound on $n_2$ similar to the
one on $n_1$ in~\eqref{eq:2938545}. Then
\begin{eqnarray*}
(2\alpha - 1)(h - 1) + 1 & = & (2(\alpha_1 + \alpha_2) - 1)(h-1) + 1 \\
& = & (2\alpha_1 - 1)(h_1 - 1) + 1 + (2\alpha_2)(h-1) \\
& = & (2\alpha_1 - 1)(h_1 - 1) + 1 + (2\alpha_2 - 1)(h-1) + (h-1) \\
& \geq & (2\alpha_1 - 1)(h_1 - 1) + 1 + (2\alpha_2 - 1)(h_2-1) + 3 \\
& \geq & n_1 + n_2 + 2 > n.
\end{eqnarray*}
Thus, we obtain strict inequality in~\eqref{eq:omega2}. This completes
the case where $G$ is disconnected.

Now suppose that $G$ is connected. Recall, that we are assuming
$\alpha \geq 2$. Thus, the connectedness of $G$, implies that $G$
contains at least one induced $3$-path $P_3 : xyz$. Using the method
of Duchet and Meyniel \cite{DuchetMeyniel82}, we construct a connected
dominating set $D \subseteq V(G)$ of $G$ with $|D| = 3 + 2k$ and
$\alpha (G[D]) \geq 2 + k$ for some non-negative integer $k$. First,
add the vertices $x,y$ and $z$ to $D$, so that $|D| = 3$ and $\alpha
(G[D]) = 2$. Obviously, $G[D]$ is connected.  If $D$ dominates $G$,
then we are done. Otherwise, since $G$ is connected, there must be
some vertex $x_1$ at distance two from $D$.  Let $z_1$ denote a vertex
adjacent to $x_1$ and a vertex of $D$, and add both $x_1$ and $z_1$ to
$D$. The induced subgraph $G[D]$ remains connected, while $|D| = 3 +
2$ and $\alpha (G[D]) \geq 2 + 1$.  Continue in this manner, at each
step adding to $D$ a vertex $x_i$ at distance two from $D$ and vertex
$z_i$ connecting $x_i$ and $D$, until we obtain a connected dominating
set $D$ with $|D| = 3 + 2k$ and $\alpha (G[D]) \geq 2 + k$ for some
non-negative integer $k$.  The process obviously terminates, since $G$
is finite.

Let $G' := G - D$.  Let $\alpha' := \alpha (G')$, $n' := n(G')$, $h' =
h(G')$ and $\alpha_D = \alpha (G[D])$. Clearly, $\alpha'\leq \alpha$
and $2 + k \leq \alpha_D \leq \alpha$. Moreover, $n' = n - 3 - 2k$,
$|D| = 3 + 2k \leq 3 + 2(\alpha - 2) = 2\alpha - 1$, since $2 + k \leq
\alpha$.

Suppose $G'$ is edge-empty. Then $\alpha' = n'$, and we establish
strict inequality in~\eqref{eq:omega2} as follows:
\begin{eqnarray*}
  (2\alpha - 1)(h - 1) + 1 & \geq & (2\alpha - 1)3 + 1 \\
  & \geq & 6\alpha - 2 \\ 
  & \geq & \alpha' + 2\alpha + 3\alpha + 1 \\
  & \geq & n' + |D| + 3 \alpha + 2 > n.
\end{eqnarray*}
Suppose that $G'$ contains at least one edge. Then, by induction,
\begin{equation}
(2\alpha' - 1)(h' - 1) + 1 \geq n'.
\label{eq:928567346}
\end{equation}
Since $G[D]$ is connected, and $D$ dominates $G$, $D$ may be
contracted to a single vertex, which dominates all other remaining
vertices, i.e., the vertices of $G'$. This observation implies $h \geq
h' + 1$. Using this fact and~\eqref{eq:928567346}, we obtain
\begin{eqnarray}
\nonumber(2 \alpha - 1) (h - 1) + 1 & \geq & (2\alpha - 1)h' + 1 \\
\nonumber         &  =     & (2 \alpha' - 1)(h' - 1) + 1 + (2\alpha - 1) \\
\label{eq:23456}  &  \geq  & n' + |D| = n.
\end{eqnarray}
Equality in~\eqref{eq:23456} implies $h=h'+1$, $\alpha' = \alpha$,
$2+k = \alpha$ and equality in~\eqref{eq:928567346}, which, by
induction, implies $h' = 2$. Thus, equality in~\eqref{eq:23456} is
impossible, since, by assumption, $h \geq 4$. This completes the
proof.
\end{pf}

\section{Proof of the Main Theorem}
\label{sec:ProofOfMainTheorem}

\begin{pf}
  Firstly, if $G$ is a graph as described in (i) or (ii), then it
  follows from Lemma~\ref{lm:lemma1} and Lemma~\ref{lm:lemma2}, that
  equality is obtained in~\eqref{eq:elem}.
  
  We prove, by induction on the order of the graph, that the
  inequality~\eqref{eq:elem} holds and that equality is attained only
  for graphs as described in (i) and (ii). The proof will be
  partitioned into several cases. Let $h := h(G)$, $\omega :=
  \omega(G)$ and $\alpha := \alpha (G)$. Obviously, the parameters $h,
  \omega $ and $\alpha $ must all be at least two, since $G$ is
  neither complete nor edge-empty. According to
  Section~\ref{sec:Heq2orHeq3}, the desired result holds for $h \in
  \{2, 3 \}$. For $\alpha = 2$ or $\omega = 2$, the desired result
  follows immediately from Theorem~\ref{th:alpha2} and
  Theorem~\ref{th:omega2}, respectively. Hence we may assume $h \geq
  4$, $\alpha \geq 3$, $\omega \geq 3$, which implies $(2\alpha -
  1)(h-1) + 3 \geq 18$, hence the desired strict inequality holds when
  $n+ \omega \leq 17$.  Thus, the base for the induction is
  established, and we shall be assuming $n+ \omega \geq 18$. We shall
  often find it convenient to introduce graphs denoted $G'$, $G_1$ and
  $G_2$; unless otherwise stated, we define $h' = h(G')$, $\omega' :=
  \omega(G')$, $\alpha' := \alpha(G')$, $n' := n(G')$ and define
  $h_i$, $\omega_i$, $\alpha_i$ and $n_i$ for $i=1,2$ analogously.

\emph{Case 1.} Suppose $G$ is disconnected, and let $G_1$ denote a
component of $G$ and let $G_2 := G - G_1$. Since $G$ is not
edge-empty, we may assume that $G_1$ is not edge-empty.
Observe that not both $G_1$ and $G_2$ are complete graphs, since that
would imply $\alpha = 2$, which contradicts our assumption $\alpha
\geq 3$.

\emph{Case 1.a.} Suppose $G_2$ is edge-empty and thus contains an
isolated vertex $x$. Define $G' := G - x$.  Obviously, $\alpha' =
\alpha - 1$, $\omega' = \omega \geq 3$ and $h' = h \geq 4$. Thus, the
induction hypothesis applies to the non-complete graph $G'$, and we obtain
\begin{eqnarray}
\nonumber (2 \alpha - 1) (h-1) + 3 & \geq & (2 + 2\alpha' - 1)(h-1) + 3 \\
\nonumber & = & 2(h-1) + (2\alpha' - 1)(h-1) + 3 \\
\nonumber & \geq & 2(h-1) + n' + \omega' \\
\label{eq:8972463} & \geq & 6 + n - 1 + \omega > n + \omega,
\end{eqnarray}
where we used the fact that $h \geq 4$. Hence we may assume that
neither $G_1$ nor $G_2$ is edge-empty.

\emph{Case 1.b.} Suppose that one of $G_1$ and $G_2$ is a complete
graph, say $G_1$. Now, by assumption, $G_2$ is neither edge-empty nor
complete, and so the induction hypothesis implies
\begin{equation}
(2 \alpha_2 - 1)(h_2 - 1) + 3 \geq n_2 + \omega_2.
\label{eq:93245783}
\end{equation}
Moreover, $\alpha = \alpha_2 + 1$, $h = \max \{ h_1, h_2 \}
\geq n_1$ and $\omega = \max \{ \omega_1, \omega_2 \}$, which allows us to establish~\eqref{eq:elem}.
\begin{eqnarray*}
(2 \alpha - 1)(h - 1) + 3 & = & (2 + 2 \alpha_2 - 1)(h - 1) + 3 \\
 & \geq & 2(h-1) + (2\alpha_2 - 1)(h_2 - 1) + 3 \\ & \geq & 2(h-1) +
 n_2 + \omega_2 \\ & = & h + n_2 + h + \omega_2 - 2 \\
 & \geq & n_1 + n_2 + h + \omega_2 - 2 \geq n + \omega,
\end{eqnarray*}
since $h \geq \omega$ and $\omega_2 \geq 2$. Observe that equality
in~\eqref{eq:elem} implies $h_2 = h = \omega$, $\omega_2 = 2$ and
equality in~\eqref{eq:93245783}. Equality in~\eqref{eq:93245783}
implies that $G_2$ is as described in Lemma~\ref{lm:lemma1} or Lemma
\ref{lm:lemma2}, and hence has $h_2 = \omega_2$. However, $h_2 =
\omega_2$ implies $h = 2$, which is impossible, since $h \geq 4$.

\emph{Case 1.c.} Suppose that neither $G_1$ nor $G_2$ is edge-empty or
complete. Now the induction hypothesis applies to both $G_1$ and
$G_2$. Obviously, $\alpha = \alpha_1 + \alpha_2$ and $\omega_1 +
\omega_2 \geq \omega + 2$. Using this and adding the
inequalities~\eqref{eq:elem} for $G_1$ and $G_2$, we obtain the
following inequality.
$$(2 \alpha - 1)(h-1) + 3 \geq n + \omega + h - 2 > n + \omega,$$
where we used the fact that $h \geq 4$.

\emph{Case 2.} Suppose that $G$ is connected and contains a complete
subgraph $K$ such that $\alpha (G-K) < \alpha (G)$. Let $G' := G - K$.
Note that $\alpha' = \alpha - 1 \geq 2$, that is, $G'$ is not a
complete graph.

\emph{Case 2.a.} Suppose $G'$ is edge-empty, that is, $\alpha' = n'$.
Then any complete subgraph of $G$ consists of some vertices of $K$ and
at most one vertex of $G'$. If, for some vertex $x \in V(G')$, the
induced graph $G[ V(K) \cup \{ x \}]$ is complete, then any
independent set of $G$ contains at most $n'$ vertices, which
contradicts the assumption $\alpha = n' + 1$. This shows that $\omega
(G) = n(K)$, and so $n(K) \geq 3$, which allows us to establish strict
inequality in~\eqref{eq:elem}.
\begin{eqnarray*}
(2 \alpha - 1)(h-1) + 3 & \geq & (2 (n' + 1) - 1) (n(K) - 1) + 3 \\
& = & (2n' + 1)(n(K) - 1) + 3 \\
& = & n' (n(K) - 1) + (n' + 1)(n(K) - 1) + 3 \\
& \geq & 2(n(K) - 1) + (n' + 1)2 + 3 \\
& = & 2n(K) + n' + n' + 3 > 2n(K) + n' = n + \omega.
\end{eqnarray*}

\emph{Case 2.b.} Suppose that $G'$ contains at least one edge, that is, $\omega' \geq 2$. Since
also $\alpha' \geq 2$, the induction hypothesis implies
\begin{equation}
(2 \alpha' - 1) \cdot (h' - 1) + 3 \geq n' + \omega',
\label{eq:92610774}
\end{equation}
which we use in the following calculations.
\begin{equation}
\begin{array}{lll}
(2 \alpha - 1) (h  - 1) + 3 & \geq & 2(h-1) + (2\alpha' - 1)(h -
1) + 3 \\
& \geq & 2(h-1) + (2\alpha' - 1)(h' - 1) + 3 \\
& \geq & 2(h-1) + n' + \omega' \\
& = & 2(h-1) + n - n(K) + \omega' \\
& = & n + h + ( \omega' - 2) + (h - n(K)) \geq n + \omega,
\end{array}
\label{eq:97237201}
\end{equation}
where equality is obtained only if we have equality
in~\eqref{eq:92610774}, $\omega' = 2$, $h = n(K)$ and $h = h' =
\omega$.  Since $h \geq 4$, we find that also $n(K) > \omega'$.
Equality in \eqref{eq:92610774} implies that $G'$ is as described in
Lemma~\ref{lm:lemma1} or Lemma \ref{lm:lemma2}, and, in particular,
$h' = \omega'$, which implies $h = 2$, contradicting the assumption $h \geq 4$.

\emph{Case 3.} Suppose that $G$ is connected and $\alpha( G - F ) =
\alpha (G)$ for every complete subgraph $F$ of $G$. Let $K$ denote a
complete subgraph of $G$ of order $\omega$, and define $G' := G - K$.
Furthermore, suppose that $G'$ is disconnected. Note that $\alpha' =
\alpha \geq 3$.

\emph{Case 3.a}
If $G'$ is edge-empty, then $n' = \alpha' = \alpha \geq 3$. Moreover,
$n = n' + n(K)$, $\omega = n(K) \geq 3$, and so we obtain strict
inequality in~\eqref{eq:elem}.
\begin{eqnarray*}
(2 \alpha - 1)(h-1) + 3 & \geq & (2 n' - 1) ( \omega - 1) + 3 \\
&   =   & (2n - \omega - 1)(\omega - 1) + 3 \\
&   =   & n(\omega - 1) + (n - \omega - 1)(\omega - 1) + 3 \\
& \geq  & 2n + (n - \omega - 1)(\omega - 1) + 3 \\
&  >   & 2n > n + \omega.
\end{eqnarray*}

\emph{Case 3.b.} Suppose that $G'$ contains at least one edge. Since
$G'$ is assumed to be disconnected, $G'$ can be partitioned into two
disjoint non-empty graphs $H_1$ and $H_2$ such that there is no edge
in $G'$ that connects a vertex of $V(H_1)$ to a vertex of
$V(H_2)$. Note that $\alpha = \alpha (G') =
\alpha (H_1) + \alpha(H_2)$.

For every vertex $x \in V(K)$, if the induced subgraphs $G [ x \cup
V(H_1) ]$ and $G [ x \cup V(H_2) ]$ have independent sets of size $a(H_1)
+ 1$ and $a(H_2) + 1$, respectively, then these independent sets can be
combined into an independent set of $G$ of cardinality $\alpha(H_1) +
\alpha(H_2) + 1$, which contradicts the fact that $\alpha = \alpha(H_1) +
\alpha(H_2)$.

Obviously, $\alpha( G [ x \cup V(H_1) ] ) \geq \alpha(H_1)$ and
$\alpha ( G [ x \cup V(H_2) ] ) \geq \alpha(H_2)$. If $\alpha( G [ x
\cup V(H_1) ] ) = \alpha(H_1)$, then we refer to $x$ as a type 1
vertex, otherwise $\alpha( G [ x \cup V(H_2) ] ) = \alpha(H_2)$, and
we refer to $x$ as a type 2 vertex. Let $T_1$ and $T_2$ denote the set
of type 1 and type 2 vertices, respectively. Then $T_1 \cup T_2 =
V(K)$ and $\alpha ( G [ T_i \cup V(H_i) ] ) = \alpha(H_i)$ for
$i=1,2$. Define $G_i := G [ T_i \cup V(H_i) ]$. Observe that $\alpha_i
= \alpha (G_i) = \alpha (H_i)$ for $i=1,2$. Clearly, the graphs $G_1$
and $G_2$ are not both edge-empty, since $H_1$ and $H_2$ are not.

\emph{Case 3.b.1.} Suppose that at least one of $G_1$ and $G_2$ is a
complete graph, say $G_1$ is a complete graph. Then $K' := G_1$ is a
complete subgraph of $G$ for which $\alpha (G-K') = \alpha_2 = \alpha
- 1 < \alpha(G)$, since $G-K' = G_2$. This contradicts the assumption
that $\alpha (G - F) = \alpha (G)$ for every complete subgraph $F$ of
$G$.

\emph{Case 3.b.2.} Now suppose that one of $G_1$ and $G_2$ is
edge-empty, while the other is neither complete nor edge-empty, say
$G_1$ is edge-empty, and $G_2$ is neither complete nor
edge-empty. This implies $T_1 = \emptyset$, $T_2 = V(K)$, $h_1 = 1$,
$\alpha_1 = n_1$ and $\omega_2 = \omega$. Moreover, the induction
hypothesis may be applied to $G_2$, that is,
$$(2\alpha_2 - 1)(h_2 - 1) + 3 \geq n_2 + \omega_2,$$
and so we obtain strict inequality in~\eqref{eq:elem}.
\begin{displaymath}
\begin{array}{rcl}
(2\alpha - 1)(h - 1) + 3 & \geq & (2(\alpha_1 + \alpha_2) - 1)(h_2 -
1) + 3 \\
& \geq & 2 \alpha_1 (h_2 - 1) + (2 \alpha_2 - 1)(h_2 - 1) + 3 \\ 
& \geq & 2 \alpha_1 (h_2 - 1) + n_2 + \omega_2 \\
& \geq & 2 n_1 (h_2 - 1) + n_2 + \omega > n + \omega
\end{array}
\end{displaymath}
since $n_1 + n_2 = n$, $n_1 \geq 1$ and $h_2 \geq \omega_2 = \omega
\geq 2$.

\emph{Case 3.b.3.} Now suppose that neither $G_1$ nor $G_2$ is
edge-empty or complete. We may assume $h_1 = h \geq 4$. Observe that
$n = n_1 + n_2$ and $\omega_1 + \omega_2 \geq |T_1| + |T_2| = \omega$.
The induction hypothesis implies
\begin{equation}
(2 \alpha_i - 1)(h_i - 1) + 3 \geq n_i + \omega_i \textrm{
  for } i = 1,2,
\label{eq:0917641}
\end{equation}
which is used in the following.
\begin{eqnarray*}
(2\alpha - 1)(h - 1) + 3 & = & (2\alpha_1 + 2\alpha_2 - 1)(h - 1) + 3 \\
 & \geq & 2\alpha_1(h - 1) + (2\alpha_2 - 1)(h_2 - 1) + 3 \\
& =     & (2\alpha_1 - 1)(h - 1) + (h - 1) + (2\alpha_2 - 1)(h_2 - 1) + 3 \\
& \geq & (2\alpha_1 - 1)(h_1 - 1) + 3 + (2\alpha_2 - 1)(h_2 - 1) + 3 \\
& \geq & n_1 + \omega_1 + n_2 + \omega_2 \geq n + \omega.
\end{eqnarray*}
Thus, equality in~\eqref{eq:elem} implies equality
in~\eqref{eq:0917641}, while both $h=4$ and $h=h_1=h_2$. Now the
induction hypothesis applied to $G_1$ and $G_2$ implies that both
$G_1$ and $G_2$ are graphs of the type described by
Lemma~\ref{lm:lemma2}. By Lemma~\ref{lm:lemma2}, $\omega_i = h_i = n_i
/ 2$ and $\alpha_i = 2$. Hence $n_i = 8$ and $\omega_i = 4$ for
$i=1,2$. Now
$$4 = h \geq \omega \geq \max \{ \omega_1, \omega_2 \} = 4$$
and so $\omega = 4$. Also, $n = n_1 + n_2 = 16$ and $\alpha = \alpha'
= \alpha_1 + \alpha_2 = 4$. Substituting these values of the
parameters into \eqref{eq:elem} we obtain the desired strict
inequality.

\emph{Case 4.} Suppose that $G$ is connected and $\alpha( G - F
) = \alpha (G)$ for every complete subgraph $F$ of $G$. Let $K$ denote
a complete subgraph of $G$ of order $\omega$, and let $G' := G - K$.
Finally, we consider the case where $G'$ is
connected.

Obviously, $G'$ is not vertex-empty. Let $x$ denote a vertex of $G'$.
According to an argument of Duchet and Meyniel~\cite{DuchetMeyniel82} (a
similar argument was given in the proof of Theorem~\ref{th:omega2}),
there exists a set $T \subseteq V(G')$ such that $x \in T \subseteq V(G')$, $|T| = 2\alpha_T - 1$, $\alpha (G[T]) = \alpha_T$, $T$ dominates $G'$ and $G[T]$ is connected.
\emph{Case 4.a.} Suppose that $T$ dominates $G$. Contract $T$, in $G$,
into one vertex $t$ and denote the resulting graph $G''$, and $J :=
G'' - t$. The vertex $t$ of $G''$ dominates $G''$, and so $h(G'') =
h(J) + 1$. Since $K \subseteq J \subseteq G$, of course $\omega(J) =
\omega$, in particular, $J$ is not edge-empty.

\emph{Case 4.a.1.} Suppose $T \subsetneq V(G')$. If $J$ were complete,
then it would contain a clique of order at least $\omega + 1$, which
is impossible, since $J \subseteq G$. Hence we may apply the induction
hypothesis and obtain the following bound.
\begin{equation}
(2\alpha (J) - 1)(h(J) - 1) + 3 \geq n(J) + \omega(J).
\label{eq:01745634}
\end{equation}
Since $h(G) \geq h(G'') \geq h(J) + 1$, we obtain
\begin{eqnarray}
\nonumber (2 \alpha - 1)(h - 1) + 3 & \geq & (2\alpha - 1)h(J) + 3 \\
\nonumber & \geq & (2\alpha(J) - 1)(h(J) - 1) + (2\alpha - 1) + 3 \\
\nonumber & \geq & n(J) + \omega(J) +  (2\alpha - 1) \\
\nonumber & =    & n(J) + \omega(J) +  |T| \\
& \geq & n + \omega,
\label{eq:9823548} 
\end{eqnarray}
where we used $\omega (J) = \omega$, $\alpha \geq \alpha_T$ and
$\alpha \geq \alpha(J)$. Thus, the desired inequality is established.
It follows from the inequalities of~\eqref{eq:9823548} that equality
in~\eqref{eq:elem} implies equality in~\eqref{eq:01745634}, $\alpha
= \alpha_T =
\alpha(J)$ and $h (G) = h(G'') = h(J) + 1$. If this is the case, then
the induction hypothesis implies that $J$ is either a graph as
described in Lemma~\ref{lm:lemma1} or Lemma~\ref{lm:lemma2}. If $J$
were a graph as described in Lemma~\ref{lm:lemma1}, then we would have
$h(J) = \omega(J) = \omega = 2$, which contradicts the assumption
$\omega \geq 3$. On the other hand, if $J$ is a graph as described in
Lemma~\ref{lm:lemma2}, then, in particular, $\alpha (J) = 2$, which
contradicts $\alpha (J) = \alpha \geq 3$. Hence, strict inequality
in~\eqref{eq:elem} is established.

\emph{Case 4.a.2.} Suppose $T = V(G')$. By contracting the vertices of
$T$ into one vertex we find that $G$ contains $K_{\omega + 1}$ as a
  minor, i.e., $h \geq \omega + 1$. Thus, in order to establish strict
  inequality in~\eqref{eq:elem} it suffices to show
  $$\frac{n + \omega - 3}{2\alpha - 1} + 1 < \omega + 1.$$
  In the
  following computations we use that fact that $\alpha = \alpha (G') =
  \alpha_T$, $|T| = 2\alpha_T - 1$ and therefore $n = \omega + 2\alpha
  - 1$.
\begin{eqnarray*}
\frac{n + \omega - 3}{2\alpha - 1} + 1 & = & \frac{2}{2\alpha - 1}
 \omega + \frac{2(2\alpha - 1) - 3}{2 \alpha - 1} \\
& < & \frac{2}{2\alpha - 1} \omega + 2 < \omega + 1,
\end{eqnarray*}
where the last strict inequality follows from the facts that $\alpha
\geq 3$ and $\omega \geq 3$.

\emph{Case 4.b.} Finally, suppose that $T$ does not dominate $G$. The
set $T$ still dominates $G'$, so $T$ does not dominate $K$. In
particular, there exists a vertex $z \in V(K)$ such that $\{ z \} \cup
S$ is an independent set of $G$ for any $\alpha_T$-set. Thus,
$\alpha_T < \alpha$. Since $G$ is connected, there exists some vertex
$x \in V(G')$ adjacent to some vertex $y \in V(K)$.  We may assume
that $T$ was construct so as to contain $x$. Since $T$ dominates $G'$,
and $y$ dominates $K$, the set $\{ y \} \cup T$ dominates all of $G$.

Contract $\{ y \} \cup T$ into one vertex $t_y$ and denote the
resulting graph $G''$, and $J := G'' - t_y$. Now $K - y \subseteq J$,
and so, since $|V(K)| \geq 3$, $J$ must contain at least one edge.

\emph{Case 4.b.1} Suppose that $J$ is not complete. Then the induction
hypothesis applies to $J$, that is,
\begin{equation}
(2\alpha(J) - 1)(h(J) - 1) + 3 \geq n(J) + \omega(J).
\label{eq:978674634}
\end{equation}
Since the vertex $t_y$ of $G''$ dominates $G''$, any minor of $J$ can
be extended to include $t_y$ in $G''$, i.e., $h(G'') \geq h(J) +
1$. Since also $h \geq h(G'')$, we obtain
\begin{eqnarray*}
(2\alpha - 1)(h - 1) + 3 & \geq & (2\alpha - 1)h(J) + 3 \\
& \geq & (2\alpha - 1)(h(J) - 1) + (2\alpha - 1) + 3 \\
& \geq & (2\alpha(J) - 1)(h(J) - 1) + (2\alpha - 1) + 3 \\
& \geq & n(J) + \omega(J) + (2\alpha - 1)  \\
& \geq & (n - |T| - 1) + (\omega - 1) +  (2\alpha - 1) \\
& \geq & (n - |T| - 1) + (\omega - 1) +  (|T| + 2) = n + \omega,
\end{eqnarray*}
where we used $\omega (J) \geq \omega - 1$, $n(J) = n - |T| - 1$,
$\alpha \geq \alpha(J)$ and $2\alpha - 1 \geq |T| + 2$. Now equality
in~\eqref{eq:elem} implies, in particular, $h = h(G'') = h(J) + 1$,
$\alpha = \alpha(J)$ and equality in~\eqref{eq:978674634}. By
induction, equality in~\eqref{eq:978674634} implies that $J$ is as
described in Lemma~\ref{lm:lemma1} or Lemma~\ref{lm:lemma2}, in
particular, $h(J) = 2$ or $\alpha(J) = 2$, which both are impossible,
since $\alpha (J) = \alpha \geq 3$ and $h(J) = h - 1 \geq 3$, by
assumption.  Thus, we obtain strict inequality in~\eqref{eq:elem}.

\emph{Case 4.b.2} Suppose $J$ is a complete graph. Recall that $\alpha
(T) < \alpha = \alpha (G')$, and so $V(G') \backslash V(T) \neq \emptyset$. Since
$V(K) \backslash \{ y \} \subseteq V(J)$ and $V(G') \backslash V(T) \subseteq V(J)$, it follows
that $V(G') \backslash V(T)$ must contain exactly one vertex, say $q$. Since
$\alpha_T < \alpha'$, it follows that any $\alpha'$-set $S$ in $G'$ must
contain the vertex $q$.  Now $S \backslash \{ q \}$ is also an independent set
in $G[T]$, and so $\alpha_T = \alpha' - 1 = \alpha - 1$.  Thus, we
obtain
$$n = n(K) + |T| + ( n' - |T| ) = \omega + (2\alpha_T - 1) + 1 = \omega + 2(\alpha - 1).$$
Now we are ready to establish the desired strict inequality.
\begin{eqnarray*}
 (2 \alpha - 1)(h - 1) + 3 & = & (2\alpha - 1)(h-2) + (2\alpha - 1) + 3 \\
&    = & (2\alpha - 1)(h-2) + (n - \omega - 1) + 3 \\
& \geq & 5(h-2) + n - \omega + 2 \\
&    = & n + h + 4h - \omega - 8 \\
& \geq & n + h + 3h - 8 > n + h \geq n + \omega,
\end{eqnarray*}
where the last strict inequality holds since $h \geq 4$. This
completes the proof.
\end{pf}

\section{Concluding Remarks}
The Duchet-Meyniel Theorem is open for further extensions and
improvements, but these will require new ideas. Some such
improvements have been obtained, for example by Wood \cite{Wood},
who proved that for $h(G)= h \geq 5$ the following inequality holds:
$$(2\alpha-1)\cdot (2h-5) \geq 2n-5.$$ The case $h = 5$ follows from
the deep result by Robertson, Seymour and Thomas \cite{RST} that any
6-chromatic graph has $K_6$ as a minor, and this is the starting
point of an induction proof.

The main problem in this area is to improve the original theorem of
Duchet and Meyniel in the case when $\alpha = 2$. This problem was
first raised by P. Seymour. The more general problem of Hadwiger's
Conjecture for $\alpha = 2$ was first pointed by W. Mader as a very
interesting special case.

\subsection*{Acknowledgment}
The result of this paper was presented by the second author at the 6th
Czek-Slovak Combinatorial Conference in Prague, July 2006. We wish to
thank the organizers of the conference for a perfect arrangement. We
also wish to thank an anonymous referee of \cite{KPT} for the original
suggestion to obtain a basic elementary extension of the Duchet-Meyniel
Theorem.

\bibliographystyle{plain}
\bibliography{DM2ndRef}
\end{document}